\documentclass[letterpaper, 9 pt, conference]{ieeeconf}

\usepackage{bm}
\usepackage{cite}
\usepackage{amsmath}
\usepackage{amssymb,amsfonts}
\usepackage[ruled,linesnumbered]{algorithm2e}
\usepackage{graphicx,color}
\usepackage{mathrsfs}
\usepackage{subfigure}
\usepackage{url}
\usepackage{color}
\usepackage{array}
\usepackage{tikz}
\usepackage{mathtools}
\usepackage{multicol}
\usepackage[colorlinks = True, linkcolor = blue, citecolor = blue]{hyperref}

\usepackage{amsthm}
\usepackage{accents}

\newtheorem{theorem}{Theorem}[section]
\newtheorem{lemma}[theorem]{Lemma}

\newtheorem*{problem-non}{Problem}
\newtheorem{remark}{Remark}

\newtheorem{example}{Example}

\newcommand{\real}{\mathbb{R}}

\newcommand{\ad}{\mathrm{ad}}

\DeclareSymbolFont{bbold}{U}{bbold}{m}{n}
\DeclareSymbolFontAlphabet{\mathbbold}{bbold}

\newcommand{\R}{\mathbb{R}}

\newcommand\oprocendsymbol{\hbox{$\square$}}
\newcommand\oprocend{\relax\ifmmode\else\unskip\hfill\fi\oprocendsymbol}

\renewcommand{\theenumi}{(\roman{enumi})}

\graphicspath{{img/}}

\def\BibTeX{{\rm B\kern-.05em{\sc i\kern-.025em b}\kern-.08em
    T\kern-.1667em\lower.7ex\hbox{E}\kern-.125emX}}
\begin{document}
\title{Data-Driven Feedback Linearization using the Koopman Generator}
  \author{Darshan Gadginmath \quad  Vishaal Krishnan \quad Fabio Pasqualetti 
                \thanks{This material is based upon work supported in part by awards ONR-N00014-19-1-2264, ARO
W911NF-20-2-0267, and  AFOSR-FA9550-19-1-0235.
                                Darshan Gadginmath \href{mailto:dgadg001@ucr.edu}{(\texttt{dgadg001@ucr.edu})} and Fabio Pasqualetti \href{mailto:fabiopas@engr.ucr.edu}{\texttt{(fabiopas@engr.ucr.edu})} are with the Department of Mechanical
                                Engineering, University of California, Riverside, Riverside, CA, 92521, USA. Vishaal Krishnan \href{mailto:vkrishnan@seas.harvard.edu}{(\texttt{vkrishnan@seas.harvard.edu})} is with the School of Engineering and Applied Sciences, Harvard University, Cambridge, MA, 02138, USA. 
                                }
                }
                
\date{}         
\maketitle

\begin{abstract} This paper contributes a theoretical framework for data-driven
 feedback linearization of nonlinear control-affine systems. We unify the
 traditional geometric perspective on feedback linearization with an
 operator-theoretic perspective involving the Koopman operator. We first show
 that if the distribution of the control vector field and its repeated Lie
 brackets with the drift vector field is involutive, then there exists an output and
 a feedback control law for which the Koopman generator is finite-dimensional and
 locally nilpotent. We use this connection to propose a data-driven algorithm `{\color{black} Koopman Generator-based Feedback Linearization (KGFL)}'
 for feedback linearization. Particularly, we use experimental data to identify
 the state transformation and control feedback from a dictionary of functions
 for which feedback linearization is achieved in a least-squares sense. {\color{black}We also propose a single-step data-driven formula which can be used to compute the linearizing transformations.} When the system is feedback linearizable and the chosen dictionary is complete, our data-driven algorithm provides the same solution as model-based feedback linearization. 
 Finally, we provide numerical examples for the data-driven algorithm and
 compare it with model-based feedback linearization. We also numerically study
 the effect of the richness of the dictionary and the size of the data set on
 the effectiveness of feedback linearization.
\end{abstract}


\section{Introduction} 
\label{sec:intro}

Nonlinear control methods rooted in model-based approaches have received
considerable attention~\cite{AI:1985}. Among these techniques, feedback
linearization has emerged as a prominent strategy, offering the implementation
of straightforward linear control methodologies to nonlinear systems. However,
a notable limitation of this approach is its demand for a comprehensive
knowledge of the system dynamics. Consequently, inadequate system
identification in the context of complex, high-dimensional cyber-physical
systems can lead to poor control performance. On the contrary, machine learning
methodologies~\cite{KJH-etal:1992,AY-FLL:1995,SS:2016} offer a robust
alternative, enabling the utilization of experimental data acquired from the
system to facilitate feedback control, even in the absence of prior knowledge
regarding the underlying system's dynamics. Nevertheless, these machine
learning methods frequently fall short of providing comprehensive insights into
both their own performance and the intricate nature of the systems they operate
on. Furthermore, the full extent of their limitations remains a subject of
ongoing investigation. The pursuit of a systematic framework for nonlinear
data-driven control remains an unresolved challenge.

Recently, significant attention has been directed towards the Koopman
operator \cite{PB-SS-SH:2021} due to its capacity to furnish a global
(infinite-dimensional) linear representation of autonomous nonlinear systems.
It was shown in \cite{MOW-IGK-CWR:2015} that the Koopman operator can be
approximated in finite dimensions with data using a dictionary of observables,
which has been a notable direction of research for nonlinear systems. However,
the commonality between the two aforementioned methodologies pertains to the
concept of complete linearization, a dimension of inquiry that has hitherto
remained unexplored in the existing literature. In this work, we bridge the gap
between the conventional technique of feedback linearization and the Koopman
operator. Furthermore, leveraging this newfound connection, we develop a
data-driven methodology capable of yielding valuable insights into the dynamics
inherent to the system. \\

\noindent \textbf{Problem setup.}
We consider a continuous-time nonlinear control-affine system, with
single input, of the form:
\begin{align} \label{eq:ctrl_affine_system}
    \dot{x} &= f(x) + g(x) u,
\end{align}
where $x \in \mathbb{X} \subseteq \real^n$ is the state, $u \in \real$ 
is the control input,
and $f,g : \mathbb{X} \rightarrow \real^n$ are the
drift and control vector fields. In the data-driven setting, 
we do not have access to the drift and control vector fields~$f, g$,
but instead have access to $N$ data samples collected from a control experiment 
on System~\eqref{eq:ctrl_affine_system}. The state and control trajectory during the experiment is $\left\{x_t,u_t\right\}$ where $t \in \mathbb{R}_{\geq0}$. 
The data collected from an experiment is represented as matrices $X, U$ as follows:
\begin{align*} 
X =\left[x_0 \ x_1 \ \ldots \ x_N \right], \quad U =\left[u_1 \ u_2 \ \ldots \ u_N \right],
\end{align*} 
   where $x_i$ and $u_i$ are the sample at the $i$-th instance of the experiment. In the experiment, the control
    $u$ is assumed to be sufficiently exciting so as to provide data of good quality~\cite{CDP-PT:2021_experiments}. For instance, $u_t$ could be sampled from a Gaussian distribution.

Our objective is to transform system~\eqref{eq:ctrl_affine_system} to a target linear system $\dot{z} = Az + Bv$, where $z$ and $v$ are transformed state and control, respectively. We propose to transform the state as $z = H(x)$ and the control as $u = \alpha(x) + \beta(x) v$. We seek to identify the transformations $H, \alpha,$ and $\beta$ using the data~$X,U$. \\



\noindent \textbf{Related work.}
 A comprehensive introduction to feedback linearization can be found in \cite
 {AI:1985}. This technique provides a systematic method to identify the necessary state and
 control transformations in the model-based case. It is crucial to note that
 these transformations are dependent on the dynamics of the system and not all
 systems allow for feedback linearization. An approximate, but still model-based,
 approach for feedback linearization was proposed in~\cite{AB-JH:1996}. These
 methods cannot be employed without a prior system identification step. Several
 works that combine learning methods for feedback linearization have been
 proposed \cite{AY-FLL:1995,SS:2016,TW-etal:2019,JU-etal:2017}. The works \cite
 {AY-FLL:1995,SS:2016} primarily use neural networks to obtain state and
 control transformations, whereas \cite{TW-etal:2019} proposes a reinforcement
 learning approach. However, these methods do not provide a clear insight into
 the control and state transformations. In \cite{LF-MM-PT:2020}, a SISO full
 state-feedback linearizable system is considered and a data-driven solution is
 proposed by approximating the system using Taylor series. An extension of the
 Willems fundamental lemma for nonlinear systems is proposed in  \cite
 {MA-etal:2022}, which is used to present a predictive control methodology with
 data. However, a systematic approach to finding the state and control
 transformations in the data-driven setting for feedback linearization has not been
 addressed in the literature. In this work, we seek to establish a data-driven
 methodology for feedback linearization which not only provides a convenient
 solution but also insight into the dynamics of the system. 

The main advantage of the Koopman operator is its ability to provide a global linear representation of a nonlinear system. However, its main drawback is its
infinite-dimensional representation for only autonomous systems. {\color{black} Recent literature has focused on finite-dimensional approximations of the Koopman
operator~\cite{MOW-IGK-CWR:2015,MH-JC:2021}. Of particular interest is the
gEDMD algorithm~\cite{SK-etal:2020} which seeks a finite-dimensional approximation of the infinitesimal generator of the Koopman operator and is based on Extended Dynamic Mode Decomposition (EDMD)~\cite{MOW-IGK-CWR:2015}.
The gEDMD algorithm~\cite{SK-etal:2020} uses a dictionary of functions to lift full-state data from an autonomous system and seeks to find a linear relation in the evolution of the lifted system.}

The works in \cite
{SLB-etal:2016,MOW-etal:2016,EK-JNK-SLB:2021,LCI-RT-MS:2022} have focused on obtaining accurate
finite-dimensional approximations of this linear operator for control. While \cite{SLB-etal:2016,MOW-etal:2016,EK-JNK-SLB:2021} have extended \cite{MOW-IGK-CWR:2015} for control, \cite{LCI-RT-MS:2022} transforms the nonlinear system as a linear parameter-varying system with the control as the variable parameter.  {\color
{black}In \cite{MK-IM:2018}, linear predictors for the control-affine nonlinear
system are considered. However, crucially, the control transformations required
for exact linearization and its connection to feedback linearization are
absent. A Luenberger observer for the system's nonlinearities is proposed using the
Koopman operator in \cite{AS:2016}. Here the control is considered
as a varying parameter, and the overall system is considered as a linear
parameter varying system. Hence, existing literature that use the Koopman operator for control have crucially missed the connection to feedback linearization. Bilinearization using the Koopman operator has also
been an area of interest \cite
{SP-SEO-CWR:2020,FN-etal:2023,DG-DAP:2017}. In \cite{SP-SEO-CWR:2020}, the
nonlinear system is approximated by interpolated bilinear systems. Then a model predictive control scheme is applied to the identified interpolated bilinear
model. Probabilistic error bounds on trajectories predicted by bilinearized
models using the Koopman operator are given in \cite{FN-etal:2023}. Conditions
for global bilinearizability using the Koopman operator are given in \cite
{DG-DAP:2017}. However, it is important to note that standard linear control techniques cannot be implemented on bilinear models. } The model-based feedback linearization approach and the modern
data-driven Koopman operator approach are both linearization techniques, yet for controlled and autonomous systems, respectively. In this paper, we focus on showing a connection between these two methods and developing a data-driven scheme for nonlinear control. \\

\noindent \textbf{Contributions.} The main contributions of this paper are as follows. We first bridge the gap between the geometric framework of feedback linearization and the Koopman operator-theoretic framework. 
In particular, we show that, when the system is involutive
to a certain degree, there exists an observable~$h$ and a feedback control~$\alpha$ such that the Koopman generator for the closed-loop system under the feedback~$\alpha$ is nilpotent at the observable~$h$. Furthermore, there
exists a finite-dimensional Koopman invariant subspace of the same dimension
as the involutive distribution for the system. This connection to the Koopman operator allows us to develop a data-driven method for feedback linearization,
by essentially casting the problem of data-driven feedback linearization as one of learning the closed-loop Koopman operator for the nonlinear control-affine system by a linearizing state/control transformation. 
{\color{black}To this end, we exploit the fact that involutivity permits a representation of the Koopman generator in the finite-dimensional Brunovsky canonical form under the linearizing state/control transformation. This allows us 
to fix the Brunovsky canonical form as the target linear representation
and learn the linearizing transformation using a set of fixed dictionary functions by a least-squares method in our algorithm \emph{Koopman generator-based Feedback Linearization (KGFL)}. We also provide a numerical feedback linearization scheme with only input-output data. With input-output data, we show that the control
 transformations can be learned in a least-squares sense using a simple
 data-driven formula. The results in \cite{CDP-DG-FP-PT:2023}, which were developed independently from this work, deal with data-driven feedback linearization with complete dictionaries for fully feedback linearizable systems. In our work, we neither make the assumption of full feedback linearizability nor of complete dictionaries. However, when the system is feedback linearizable and the dictionaries used in KGFL are complete, the solution is exact and is equal to the model-based solution.} Finally, we demonstrate the performance of our algorithm
 with numerical simulations for multiple examples. We perform both full state
 feedback linearization and output feedback linearization on the Van der Pol 
 oscillator {\color{black} and compare it against existing nonlinear data-driven
 control techniques. We consider a higher dimensional system with the control
 entering nonlinearly and show that our algorithm can be used for complex
 systems.} We also provide insight on the effect of richness of dicitonary and
 data size on the accuracy of feedback linearization method. 

\section{Preliminaries}
Let $\real^d$ be the $d$-dimensional Euclidean space.
Let $(\mathbb{X}, d_{\mathbb{X}})$ and $(\mathbb{Y}, d_{\mathbb{Y}})$ be metric spaces. A map $f : \mathbb{X} \rightarrow \mathbb{Y}$ is said to be \emph{Lipschitz} (with Lipschitz constant $\ell_f$) if $d_{\mathbb{Y}}(f(x_1), f(x_2)) \leq \ell_f d_{\mathbb{X}}(x_1, x_2)$ for any $x_1, x_2 \in {\mathbb{X}}$.
The space of $k$-times continuously differentiable functions on $\mathbb{X}$ is denoted by $C^k(\mathbb{X})$.
Let $V$ be a normed vector space and let $T:V\rightarrow V$
be a bounded linear operator on~$V$. A subspace $W \subseteq V$ 
is said to be \emph{$T$-invariant} if $T(W) \subseteq W$. 
The operator $T: V\rightarrow V$ is \emph{locally nilpotent} 
with index $r$ at $v \in V$ if $T^k v \neq 0$ for all 
$k \in \lbrace 0, \ldots, r-1 \rbrace$ and $T^r v = 0$.
Furthermore, if $v, T(v), \ldots, T^r(v)$ are linearly 
independent, then $\mathrm{span}\{v,T(v),\ldots,T^r(v)\}$ is 
said to be a \emph{$T$-cyclic subspace} of~$T$.
For a measure space $(\mathbb{X}, \Sigma, \mu)$ (where $\Sigma$ is a sigma-algebra on $X$ and $\mu$ is the measure on $(X, \Sigma)$), a property $P$ is said to hold \emph{almost everywhere (a.e.)} if the subset over which the property $P$ fails to hold is of $\mu$-measure zero. The \emph{Lie bracket} between two vector fields $f$ and $g$ 
is denoted by $[f,g]= \ad_f g = L_g f - L_f g$. The \emph{adjoint} of order $k$ is defined recursively as 
$\ad_f^k g = \left[ f \; , \; \ad_f^{(k-1)} g \right]$ with $\ad_f^0 g = g$. {\color{black}The Gateaux derivative $U(v; \eta)$~\cite{WR:64} of operator $T \in C^1(V,V)$ at $v \in V$ along $\eta \in V$
is given by
\begin{align*}
    \lim_{h \rightarrow 0^+} \left\| T(v+h\eta) - T(v) - U(v; \eta) h \right\|_{V} = 0.
\end{align*}
$W \subseteq V$ is called a \emph{stable subspace} of~$T$ with respect to  perturbations along $\eta$ if $U(w; \eta) = 0$ for all $w \in W$.}

\subsection{Lie derivative as Koopman generator}
Consider the autonomous system $\dot{x}(t) = f(x(t))$ with state space $\mathbb{X} \subset \real^d$,
where $f: \mathbb{X} \rightarrow \real^d$ is Lipschitz.
Let $\mathbf{s} : \mathbb{X} \times \real_{\geq 0} \rightarrow \mathbb{X}$ be the flow of the vector field $f$, such that for any $x \in \mathbb{X}$, $\mathbf{s}(x,0) = x$ and $\frac{d}{dt} \mathbf{s}(x, t) = f \left( \mathbf{s}(x,t) \right)$. 
For a function $\phi \in C^1(\mathbb{X})$, let $L_f \phi$
be the Lie derivative of $\phi$ with respect to~$f$, such that
for any $x \in \mathbb{X}$, $L_f$ uniquely satisfies
\begin{align*}
    \lim_{h \rightarrow 0^+} \left| \phi(\mathbf{s}(x,h)) - \phi(x)  - h (L_f \phi)(x) \right| = 0.
\end{align*}
Let $\mathcal{K}: C^1(\mathbb{X}) \times \real_{\geq 0} \rightarrow C^1(\mathbb{X})$ 
be the Koopman operator for the flow $\mathbf{s}$, such that for any $x \in \mathbb{X}$ we have:
\begin{align*}
    (\mathcal{K}_t \phi) (x) = \phi(\mathbf{s}(x,t)),
\end{align*}
where we have adopted the notation $\mathcal{K}_t \phi = \mathcal{K}(\phi, t)$.
The family of operators $\lbrace \mathcal{K}_t \rbrace_{t \geq 0}$ 
is said to be right-differentiable at $t=0$ with derivative operator $\mathcal{L}$
if the following holds for any $\phi \in C^1(\mathbb{X})$:
\begin{align*}
    \lim_{h \rightarrow 0^+} \left\| \mathcal{K}_h \phi - \mathcal{K}_0 \phi - h \mathcal{L} \phi \right\|_{L^2(\mathbb{X})} = 0.
\end{align*}
In operator-theoretic terminology, 
$\mathcal{L}$ is called the infinitesimal generator of 
the family $\lbrace \mathcal{K}_h \rbrace_{h \geq 0}$,
and since $\mathcal{K}$ is the Koopman operator,
we thereby refer to $\mathcal{L}$ as the \emph{Koopman generator}.
The following lemma establishes the relationship between
the Koopman generator corresponding to the flow~$\mathbf{s}$
and the Lie derivative~$L_f$ corresponding to the vector field~$f$:
\begin{lemma}[\bf \emph{Lie derivative as Koopman generator}] 
\label{lemma:Koopman-lieder} 
Let $f: \mathbb{X} \rightarrow \real^n$ be Lipschitz,
with flow $\mathbf{s}: \mathbb{X} \times \real_{\geq 0} \rightarrow \mathbb{X}$,
and let $\phi \in C^1(\mathbb{X})$. 
Then, almost everywhere in $\mathbb{X}$, the 
Lie derivative $L_f \phi$ is equal to the infinitesimal generator~$\mathcal{L}$ of the 
Koopman operator~$\mathcal{K}$ for the flow~$\mathbf{s}$, i.e., 
$L_f \phi = \mathcal{L}\phi$ a.e. in~$\mathbb{X}$~\cite[Section 7.6]{AL-MCM:1994}.
\end{lemma}
\noindent We do not include the proof of Lemma~\ref{lemma:Koopman-lieder}
for the sake of brevity. We note that Lemma~\ref{lemma:Koopman-lieder}
establishes that the Lie derivative operator and the 
Koopman generator are equivalent (in an a.e. sense) 
over the space of continuously differentiable observable functions.
{\color{black}We denote the closed-loop Koopman generator for System~\eqref{eq:ctrl_affine_system} 
with feedback $u = \alpha(x)$ as 
$L_{f + g \alpha}$. Note that the closed-loop Koopman generator is linear over the space of 
continuously differentiable functions $C^1(\mathbb{X})$, i.e. $L_{f+g\alpha} = L_f + \alpha L_g$\footnote{{\color{black}For any Lipschitz
$f_1, f_2:\mathbb{X}\rightarrow \real^n$ 
and $k_1, k_2 \in C^1(\mathbb{X})$, 
the Lie derivative of $\phi \in C^1(\mathbb{X})$ with respect to 
$k_1 f_1 + k_2 f_2$ 
is given by $L_{k_1 f_1 + k_2 f_2} \phi
= \nabla \phi \cdot \left( k_1 f_1 + k_2 f_2 \right)
= k_1 \nabla \phi \cdot f_1 + k_2 \nabla \phi \cdot f_2
= k_1 L_{f_1} \phi + k_2 L_{f_2} \phi$.}}. }

\subsection{Feedback Linearization}
Feedback linearization addresses the problem of designing a linearizing feedback controller in the model-based setting (given vector fields $f$ and $g$).
{\color{black}
 An output function $h \in C^1(\mathbb{X})$
 of system~\eqref{eq:ctrl_affine_system} is
 said to have \emph{relative degree}~$r$~\cite{AI:1985} if
 $L_gL_f^k h(x) = 0$ for $k = \{ 0,\dots, r-2 \}$ and
 $L_gL_f^{r-1} h(x) \neq 0$.}
 For an output function~$h$ of relative degree~$r$,
 the feedback linearizing transformation is given~by:
\begin{align}
z &= H(x) =  \left[ \begin{matrix}h(x) \\ L_f h(x) \\ \vdots \\ L^{(r-1)}_f h(x)\end{matrix}\right], \quad u = \alpha(x) + \beta(x) v, \label{eqn:FL-transform}
\end{align} 
yielding an $r$-dimensional linear system 
in the transformed state and control $(z,v)$ of the form,
\begin{align} \label{eqn:lin-sys-FL}
\begin{aligned}
\dot{z} = Az + Bv,  \quad \text{s.t.}~
A = \left[
\begin{matrix}
0 & 1 & 0 &\dots&0 \\ 
0 & 0 & 1 &\dots&0 \\
\vdots & \vdots & \vdots & \ddots & 0 \\
0 & 0 & 0 & \dots & 1 \\
0 & 0 & 0 & 0 & 0 
\end{matrix}\right],~ B = \left[\begin{matrix}
0 \\ 0 \\ \vdots \\ 1
\end{matrix}\right]
\end{aligned}
\end{align}
The above transformation results in an $(n-r)$-dimensional residual dynamics,  
called \emph{zero dynamics}, which is unobservable and uncontrollable,
and the system is said to be \emph{input-output feedback linearizable}. 
Therefore, this approach to linearization-based control
relies crucially on the choice of the output function $h$ 
which results in a stable zero dynamics.
If the relative degree $r=n$, the state space dimension, 
the system is said to be \emph{full-state feedback linearizable},
which is the case if and only if it is both controllable\footnote{The system~\eqref{eq:ctrl_affine_system} is said to be \emph{controllable} when the distribution $\Delta(x) = \mathrm{span} \left\{g(x), \ad_f g(x), \dots,  \ad_f^{k-1} g(x), \ldots \right\}$ is such that $\mathrm{dim}\left(\Delta(x)\right) = n$, for all $x \in \mathbb{R}^n$.} and 
integrable\footnote{A set of linearly independent vector fields $\{f_1(x),f_2(x),\dots,f_m(x)\}$ on $\real^n$ is said to be \emph{integrable}~\cite[Definition 6.4]{JJES-WL-etal:1991} if and only if there exist $n-m$ scalar functions $\{h_1(x),h_2(x),\dots,h_{n-m}(x)\}$, with $\nabla h_i$ linearly independent, satisfying $\nabla h_i \cdot f_j = 0$ for any $i \in \{ 1, \dots, n-m \}$ and $j \in \{ 1, \ldots, m \}$.}~\cite[Theorem~6.2]{JJES-WL-etal:1991}. 
Furthermore, by the Frobenius theorem~\cite{JJES-WL-etal:1991}, a distribution of 
linearly independent vector fields $f_1(x),f_2(x),\dots,f_m(x)$ is completely integrable if and only if it is involutive\footnote{The distribution $\Delta = \mathrm{span}\{f_1,f_2,\dots,f_m\}$  is said to be \emph{involutive}~\cite[Definition 6.5]{JJES-WL-etal:1991} if and only if $\ad_{f_i} f_j \in \Delta$ for any $i,j \in \{1,\dots,m\}$.}.

{\color{black}
\begin{example}[\bf \emph {Input-Output Feedback Linearizable System}]  Consider a system with the following vector fields, 
\begin{align*}
f(x) = \left[\begin{matrix} 2 x^2_3 \\ -1 \\ 0 \end{matrix}\right], \ g(x) = \left[\begin{matrix} -x_1 \\ -2x_2 \\ \frac{1}{2}x_3 \end{matrix}\right].
\end{align*}
The distribution $\Delta$ is
\begin{align*}
\Delta &= \{g,\mathrm{ad}_fg\}  = \left(\begin{matrix}-x_1 & -4x_3^2 \\ -2x_2 & 2 \\ \frac{1}{2}x_3 & 0 \end{matrix}\right).
\end{align*}
It can be noticed that $\Delta$ has a rank equal to 2 for all $x$. Further, $\ad_f^k g= 0$ for all $k \geq 2$. Thus $\{g,\ad_f g, \ad_f^k g \} $ has rank 2 for all $k$ and $x$, and $\Delta$ is involutive. However, note that the system is not controllable. Therefore, the system is not full-state feedback linearizable and only input-output feedback linearizable.
\end{example}
}

\section{Data-driven feedback linearization}
In this section, we propose a data-driven technique, called Koopman generator-based Feedback Linearization (KGFL), to perform data-driven feedback linearization to stabilize System~\eqref{eq:ctrl_affine_system}. We first establish that there exists an observable and feedback control that render the closed-loop Koopman generator finite-dimensional. We then seek to find this transformation using experimental data. 
\subsection{Koopman generator-based feedback linearization}

We now establish the connection between feedback linearization and the closed-loop Koopman generator $L_{f+g\alpha}$, which will serve as the basis for the numerical algorithm to determine the linearizing state/control transformation. To this end, the following theorem establishes the relationship between the geometry of the system, as manifested in the involutivity of the distribution $\Delta=\mathrm{span}\{g,\ad_f g,\dots,\ad_f^{r-2} g\}$, 
and a Koopman invariant subspace: 

\begin{theorem}[\bf \emph{Involutivity of distribution and nilpotency of Koopman generator}] \label{thm:Koopman-nilpotent}
Let $\Delta=\mathrm{span}\{g,ad_f g,\dots,ad_{f}^{r-2} g\}$ be an 
$(r-1)$-dimensional distribution for $f, g$ in System~\eqref{eq:ctrl_affine_system}.
The following are equivalent:
\begin{enumerate}
    \item $\Delta$ is involutive. 
    \item There exists $h \in C^r \left( \mathbb{X} \right)$ 
    and $\alpha \in C(\mathbb{X})$ such that $L_{f + g \alpha}$ is 
    locally nilpotent (with index $r$) at $h$ and the associated cyclic subspace is a stable subspace with respect to $\alpha$.
    \oprocendsymbol
\end{enumerate}
\end{theorem}
\noindent We refer the reader to Appendix~\ref{app:proof_koopman_nilpotent} for the proof of Theorem~\ref{thm:Koopman-nilpotent}. Some comments on Theorem~\ref{thm:Koopman-nilpotent} are in order.
We first note that $h$ is an observable and 
$L_{f+g\alpha}$ is the closed-loop Koopman generator with feedback $\alpha$.
From Theorem~\ref{thm:Koopman-nilpotent}, we can see the existence of a subspace invariant to the closed-loop Koopman generator $L_{f+g\alpha}$. This Koopman invariant subspace is indeed $\mathrm{span}\{h(x), \ L_{f} h(x),\ \dots ,\ L^{r-1}_{f} h
 (x) \}$. As the functions $L_f^i h(x),$ \ for all $i = \{0,1,\dots,r-1\}$, are
 linearly independent, the dimension of the Koopman invariant subspace is $r$ and the closed-loop Koopman generator $L_{f+g\alpha}$ has a finite-dimensional representation. We now consider $H(x) = \left[h(x) \ L_f h(x) \ \dots \ L^{r-1}_f h(x) \right]^\top$, whose dynamics under feedback~$\alpha$ is given by
 \begin{align*}
\frac{\mathrm{d}H(x)}{\mathrm{d}t}&= L_{f+g\alpha}H = A H(x),
 \end{align*}
where $A$ is as defined in equation~\eqref{eqn:lin-sys-FL}. 
For a system that is full-state feedback linearizable, the closed-loop Koopman generator has an $n$-dimensional invariant subspace, and Theorem~\ref{thm:Koopman-nilpotent} can be seen as a choice of an observable $h$ and feedback $\alpha$, that induces the linear system~\eqref{eqn:lin-sys-FL} of dimension $n$. 
Furthermore, Theorem~\ref{thm:Koopman-nilpotent} emphasizes the role of involutivity in ensuring the existence of a feedback $\alpha$ such that there exists a finite-dimensional invariant subspace for the closed-loop Koopman generator $L_{f + g \alpha}$.  For the special case when the open-loop Koopman generator $L_f$ already has a finite-dimensional invariant subspace, then we can simply let $\alpha = 0$. For instance, when the vector field $f$ is a constant matrix in $\mathbb{R}^{n\times n}$, $L_f$ has an invariant subspace. The role of $\alpha$ is to induce a finite-dimensional invariant subspace for $L_{f + g \alpha}$. 




{\color{black}
\begin{remark}[\bf \emph{Multiple Inputs}] When there are multiple inputs, the
 observable $h$ has a relative $r_i$ for each input $i$. If there are $m$ inputs, the
 feedback linearizing inputs will now be $\alpha \in \mathbb{R}^m$ and
 $\beta \in \mathbb{R}^{m \times m}$. Particularly, the transformed system will
 look as follows:
\begin{align*}
\dot{z} &= A_c z + B_c v
\end{align*}
where $A_c = \mathrm{blkdiag}\{A_1,A_2,\dots,A_m\}$ and $B_c = \mathrm{blkdiag}\{B_1,B_2,\dots,B_m\}$. Here $\mathrm{blkdiag}(\cdot)$ returns the block matrix with its arguments on the diagonal and block zero matrices everywhere else.
\begin{align*}
A_i &= \left[\begin{matrix}  0 & 1 & 0 & \dots & 0 \\
                             0 & 0 & 1 & \dots & 0  \\
                             \vdots & \vdots &\vdots &\ddots & \vdots \\
                             0 & 0 & 0 & \dots & 1 \\ 
                             0 & 0 & 0 & \dots & 0  
                \end{matrix}\right] \in \mathbb{R}^{r_i \times r_i} , \ B_i = \left[\begin{matrix}  0 \\
                             0 \\
                             \vdots \\ 
                             0 \\
                             1  
                \end{matrix}\right] \in \mathbb{R}^{r_i \times 1}.
\end{align*}
Here, $r_i$ is the relative degree of the observable $h$ with respect to the $i$-th input such that $r_1 + r_2 + \dots + r_m = n$. This also means that there exists a Koopman-invariant subspace with respect to every input, and also all the inputs considered simultaneously. 
\end{remark}
}

\subsection{Data-driven algorithm} \label{subsect:fl}
We use a dictionary of functions to identify the Koopman generator to lift full-state data from the nonlinear control-affine system and seek to find a linear relation in the evolution of the lifted system. Furthermore, since we already know the form of the Koopman generator, particularly the matrices~$A$ and~$B$ as defined in Equation~\eqref{eqn:lin-sys-FL} which define the action of the closed-loop Koopman generator, we seek to best approximate this structure. 
To this end, we first recall the state and control transformations~\eqref{eqn:FL-transform}, i.e.,
$z = H(x)$ and $u = \alpha(x) + \beta(x) v$.
We note that the state transformation $H$ is determined by the 
observable $h$ for which we seek an estimate $\hat{h}$,
expressed using a dictionary $\phi$ (consisting $M$ real functions, i.e., $\phi = \left[ \phi_1 \ \phi_2 \ \dots \ \phi_M\right]^\top$, where
 $\phi_i:\R^n\rightarrow \R$) 
 as $\hat{h}(x) = K^\top \phi(x)$, where $K \in \R^M$.
We then invert the control transformation to express the external control $v$ in terms of $u$ as
\begin{align*}
    v = \zeta(x) + \eta(x) u,
\end{align*}
where $\zeta$, $\eta$ are given by
\begin{align*}
    \zeta(x) = - \frac{\alpha(x)}{\beta(x)}, \quad \eta(x) = \frac{1}{\beta(x)}.
\end{align*}
Note that $\alpha, \beta$ can in turn be obtained from $\zeta, \eta$ as
$\alpha(x) = - \zeta(x)/\eta(x)$ and $\beta(x) = 1/\eta(x)$.
We now seek estimates $\hat{\zeta}$ and $\hat{\eta}$ for the above, 
which are expressed using dictionaries $\theta, \gamma$ (consisting $k$ real functions) as 
\begin{align*}
\hat{\zeta}(x) &= G^\top \theta(x), \quad \hat{\eta}(x) = J^\top \gamma(x).   
\end{align*}
where~$G, J \in \R^k$ are to be estimated from data. 
The dictionary of repeated time derivatives of the 
observable dictionary~$\phi$ is represented by D as
\begin{align*}
D(x) &= \left[\begin{matrix}
\phi(x)^\top \\
\dot{\phi}(x)^\top\\
\vdots \\ 
\phi^{(r-1)}(x)^\top
\end{matrix}\right].
\end{align*}
Here $D(x)$ utilizes the structure of $H(x)$ as the state transformation $H(x)$ contains the observable $h$ and its repeated derivatives. Utilization of this structure for the dictionary is also a novelty of our algorithm.
The estimated state transformation $\hat{z}(x)$ and control transformations can now be represented as follows:
\begin{align*}
\begin{aligned}
z = D(x) K, \qquad
v = G^\top\theta(x) + J^\top\gamma(x) u. 
\end{aligned}
\end{align*} 
We know from feedback linearization that the state and control
 transformation yield the linear system~\eqref{eqn:lin-sys-FL}. Here we emphasize the fact that we do not assume that the system is feedback linearizable or the dictionaries contain all the necessary functions for the required transformations. Therefore, we
 obtain the vectors $K,G$ and $J$ that approximate the structure of the Koopman generator in a least-squares sense, as follows:
\begin{align}
\label{prob:optim-main}
\begin{split}
&\underset{K,G,J}{\min}~\sum_{t=1}^{N-1} \left\| 
\frac{D(x_{t+1}) - D(x_t)}{\tau} K  - A \  D(x_t) K  - B v_t \right\|^2, \\
& \mathrm{s.t.} \quad
v_t = G^\top \theta(x_t) + J^\top \gamma(x_t) u_t,  \quad \forall ~t \in \{1,\dots,N\}.
\end{split}
\end{align}
\noindent Problem~\eqref{prob:optim-main} seeks to obtain the 
vectors $K, \ G$ and $J$ simultaneously.
It is evident that $K=0$, $G=0$ and an arbitrary non-zero $J$ minimizes the cost in \eqref{prob:optim-main}. Therefore, simple constraints to make $K\neq 0$ and $G\neq 0$ can be imposed. We note that Problem~\eqref{prob:optim-main} involves 
a finite difference approximation with sampling interval $\tau$ of the time derivative of $D$.

{\color{black} 
In the preceding analysis, we had assumed that the system state is directly observable. However, when we only have access to an output $y = h(x)$,
through an observable~$h$, whereby the data consists of inputs and outputs in the following matrices
\begin{align*} 
Y =\left[y_0 \ y_1 \ \ldots \ y_N \right], \quad U =\left[u_1 \ u_2 \ \ldots \ u_N \right],
\end{align*} 
the problem becomes one of input-output feedback linearization.
The problem of data-driven input-output feedback linearization is to find the necessary transformations $H,$ $\zeta$ and $\eta$ using only input-output data $Y, U$. However, the sub-problem of finding the state transformation $H(x)$ in the input-output feedback linearization problem is simpler as $H(x)$ is computed directly from data, as $y_t$ and its repeated time derivatives, since the observable~$h$ is apriori fixed,
\begin{align*}
z = H(x) = \left[y \ \dot{y} \ \dots \ y^{r-1}\right]^\top.
\end{align*}
Therefore, the data-driven input-output feedback linearization problem reduces to finding estimates $\zeta$ and $\eta$ using input-output data,
formulated as follows:
\begin{align}
\label{prob:optim-iofl}
\begin{aligned}
&\underset{G,J}{\min} ~\sum_{t=1}^N \left\| 
\frac{\mathrm{d}z_t}{\mathrm{d}t}  - A \  z_t - B v_t  \right\|^2_2, \\
&\mathrm{where} ~ \begin{cases} v_t = G^\top \theta(z_t) + J^\top \gamma(z_t) u_t,  \quad \forall t \in \{ 1, \dots , N \}  \\
z_t = \left[y_t \ \dot{y}_t \ \dots \ y^{(r-1)}_t \right]^\top
\end{cases}
\end{aligned}
\end{align}

%
We now present two methods to obtain the solutions to problems~\eqref{prob:optim-main} and~\eqref{prob:optim-iofl}.
\subsubsection{Iterative Algorithm - KGFL}
The KGFL algorithm is an iterative algorithm based on gradient descent. The gradients of the cost in Problem~\eqref{prob:optim-main} with respect to the parameters $K$, $G$ and $J$ can be explicitly computed. These gradients are utilized in an iterative manner to obtain the estimates $\hat{z}$, $\hat{\zeta}$ and $\hat{\eta}$. 
Algorithm~\ref{algo:datadriven} outlines the 
Koopman Generator-based Feedback Linearization (KGFL) algorithm.
{\color{black} The parameters $K$, $G$ and $J$ at iteration $i$ are denoted by $K(i)$, $G(i)$ and $J(i)$ respectively. In the algorithm, the state transformation parameter $K$ is computed while keeping the control transformation parameters $G$ and $J$ fixed from the previous iterations. Subsequently, the control transformation parameters $G$ and $J$ are fixed while keeping the state transformation parameters fixed. In the algorithm, we denote the cost function in Problem~\eqref{prob:optim-main} as $C$ and its gradient with a parameter $Q$ as $\nabla_Q$ C. We emphasize that the dictionaries $D$, $\theta$ and $\gamma$ may not contain all the nonlinearities of the system's dynamics. Hence, we solve the data-driven feedback linearization problem \eqref{prob:optim-main} in a least-squares sense.}  

\begin{algorithm} \label{algo:datadriven}
\caption{Koopman generator-based Feedback Linearization (KGFL)}
\textbf{Data:} $X$ and $U$ from System~\eqref{eq:ctrl_affine_system}\\
\textbf{Initialize:} Dictionaries $\phi,\theta$, $\gamma$; Number of iterations $E$; Initial guess for $K, G$ and $J$; Learning rate $\epsilon$\\
\SetAlgoLined
For $i$ = 1 to $E$:  \\
\ \ \ $K(i) =  K(i-1) -\epsilon \ \nabla_K C(K(i-1),G(i-1),J(i-1)) $ \\ 
\ \ \ $G(i) =  G(i-1) -\epsilon \ \nabla_G C(K(i),G(i-1),J(i-1)) $ \\ \label{solveGJ}
\ \ \ $J(i) =  J(i-1) -\epsilon \ \nabla_J C(K(i),G(i),J(i-1)) $
\end{algorithm}
For input-output feedback linearization, step 4 in KGFL need not be performed as the state transformation is fixed apriori.

\subsubsection{Single-step method}
The solutions to Problems~\eqref{prob:optim-main} and~\eqref{prob:optim-iofl} can also be computed in a single step. We characterize the solutions to the two problems in the following theorem:
}
\begin{theorem}[\bf \emph{Full-state 
and input-output linearizing transformations}] 
\label{thm:state_io-fl}
(a) Full-state feedback linearization:
The solution to Problem~\eqref{prob:optim-main} is given by,
\begin{align}
\left[G^\top \ J^\top \right] = B^\top \left( \frac{dD(X)}{dt} - A D(X) \right) K \left[\begin{matrix} \Theta(Y) \\ \Gamma(Y) \otimes U\end{matrix}\right]^\dagger,\label{eqn:iofl-LSsolution}
\end{align}
(b) Input-Output linearization: 
The solution to problem~\eqref{prob:optim-iofl} is given by,
\begin{align}
\left[G^\top \ J^\top \right] = B^\top(\dot{\mathbf{Z}} - A \mathbf{Z}) \left[\begin{matrix} \Theta(Y) \\ \Gamma(Y) \otimes U\end{matrix}\right]^\dagger,\label{eqn:iofl-LSsolution}
\end{align}
where, 
\begin{align*}
\mathbf{Z} &= \left[\begin{matrix}z_1 & z_2 & \dots & z_N
 \end{matrix} \right] \\
\dot{\mathbf{Z}} &= \left[\begin{matrix}\dot{z}_1 & \dot{z}_2 & \dots & \dot{z}_N
 \end{matrix} \right], \\
\Theta(Y) &= \left[\theta(z_1) \ \theta(z_2) \ \dots \ \theta(z_N)\right]  \\
\Gamma(Y) &= \left[\gamma(z_1) \ \gamma(z_2) \ \dots \ \gamma(z_N)\right] . 
\end{align*}\hfill \oprocendsymbol
\end{theorem}
\noindent We refer the reader to Appendix~\ref{app:proof_lemma_io-full}
for the proof of Theorem~\ref{thm:state_io-fl}.
Some comments on Theorem~\ref{thm:state_io-fl} are now in order. For full-state feedback linearization, \cite{CDP-DG-FP-PT:2023} presents a similar solution when the dictionaries are complete. The solution to $K$, $G$ and $J$ are represented as vectors belonging to the null space of a matrix only dependent on the data $X,U$. Further, the input-output feedback linearization case is unaddressed in \cite{CDP-DG-FP-PT:2023}. The closed-form least-squares solution provides a numerically simple method to compute the control
transformations. It is important to note that the relative degree affects the
number of repeated derivatives that are required to compute the solution
in \eqref{eqn:iofl-LSsolution}. If the relative degree of the output is not
known, the control designer can solve with a different number of repeated
derivatives and choose one that results in the least error. 


\section{Numerical Results}In this section, we demonstrate the effectiveness
 of the proposed data-driven feedback linearization technique KGFL. We run the numerical experiments on an i9-9900K CPU with 128GB of RAM. We sample $u_i$ from $\mathcal{N}(0,5)$ to collect data, and the state is sampled every
 $0.01$ seconds. We choose the control task of stabilization. The exogenous input
 to the linearized system $v$ is chosen to be $v_t = \left
 [-2 \ -2 \right]^\top z_t$ for full state feedback linearization. This places
 the poles of system~\eqref{eqn:lin-sys-FL} at $-1 \pm 1i$. For input-output
 feedback linearization with relative degree 1, we choose $v_t = -2 z_t$. Since
 the system~\eqref{eqn:lin-sys-FL} is controllable, we can arbitrarily place the poles of the system using state feedback.

The dictionaries $\phi,\theta$ and $\gamma$
 are chosen as proposed in \cite{LS-KK:2021}. Particularly, these dictionaries
 contain Kronecker products of the Hermite polynomials of individual states.
 The probabilist's Hermite polynomial $H_n(x)$ \cite{EC-MG-MADO:2021} of order $n$ is defined as $H_n
 (x) = (-1)^n e^{\frac{x^2}{2}} \frac{\mathrm{d}^n}{\mathrm{d}x^n} e^{-\frac{x^2}{2}}$. 
 For example, the Hermite polynomials for $n = \{0,\dots,4\}$ are $1,x,x^2-1,x^3-3x,x^4-6x^2+3$, respectively. Hermite polynomials serve as a useful choice for dictionaries as they form an
 orthogonal basis of the Hilbert space of functions \cite{EC-MG-MADO:2021}. 

\subsection{Numerical testbeds}\label{sec:vanderpol}
We consider two testbed systems for our numerical experiments, 
to demonstrate the proposed data-driven feedback linearization algorithm.
The first testbed system we consider is the classical Van der Pol oscillator {\color{black}but with the input entering the system nonlinearly as}
\begin{align}
\begin{split}
\dot{x}_{1} &= x_{2}, \\
\dot{x}_{2} &= -x_{1} + 0.5(1-x_{1}^2)x_{2} + (1-x_2^2)u_t.
\end{split}\label{eqn:vanderpol}
\end{align} 
We demonstrate both full-state feedback linearization and input-output feedback linearization for the Van der Pol oscillator~\eqref{eqn:vanderpol}. 
{\color{black}The second testbed system we consider is an arbitrary feedback linearizable system of 6 dimensions. The system is represented in state space as follows:
\begin{align}
\begin{aligned}
\dot{x}_{1} &= a_1 x_{2}, \ \ \dot{x}_{2} = a_2 x_{3}, \ \ \dot{x}_3 = a_3 x_4,  \\
\dot{x}_4 &= a_4 x_5, \ \ \dot{x}_5 = a_5 x_6 \\
\dot{x}_6 &= a_6 x_{1}^2 - a_7 \sin(x_{2})  + a_8 x_3^2 + a_9 x_4 ^ 2 + a_{10} (1-x_1^2) u. 
\end{aligned}\label{eqn:complex-sys}
\end{align}
where the coefficients $a_i$, with $i \in \{1,2, \dots, 10 \}$, are sampled from a standard Gaussian $\mathcal{N}(0,1)$. We demonstrate full-state feedback linearization with stabilization at the origin. For System~\eqref{eqn:complex-sys}, we choose a dictionary similar to that chosen for the Van der Pol oscillator in Section~\ref{sec:vanderpol}. Further, we augment the dictionary with $\sin(x_{1}), \sin(x_{2})$ and $\sin(x_{3})$.}

\subsubsection{Full state feedback linearization}
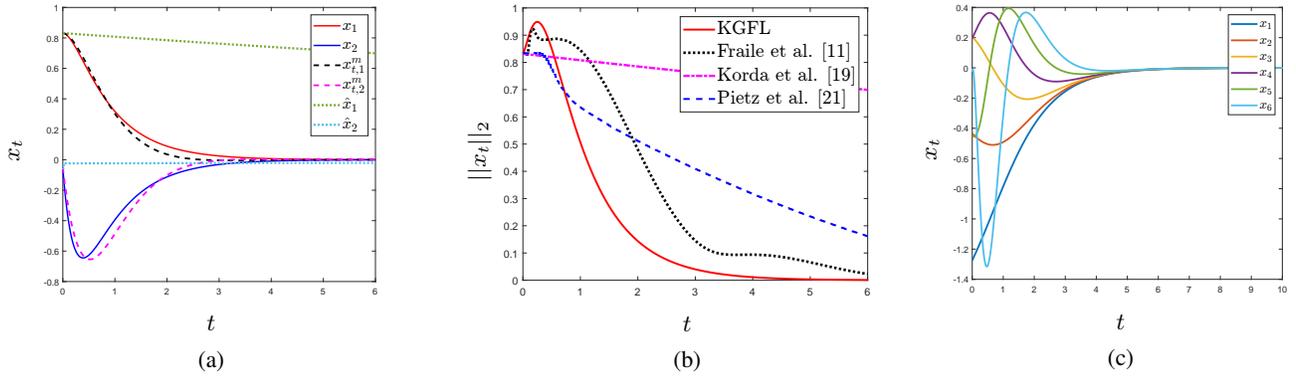
\begin{figure*}[htb]
\begin{multicols}{3}
\begin{tikzpicture}
  \node (img1)  {\includegraphics[width=0.25\textwidth]{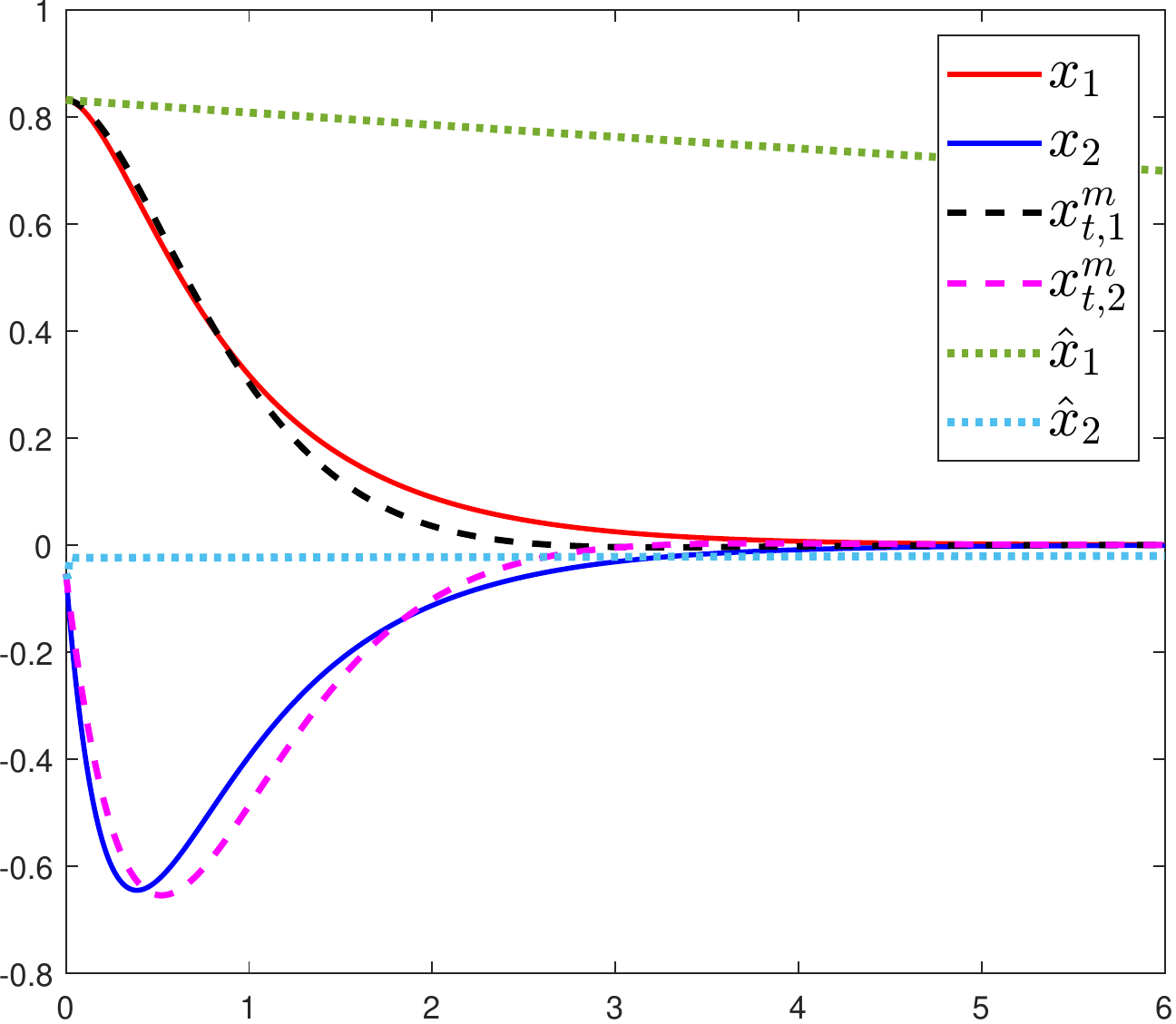}};
  \node[below of= img1, node distance=0cm, yshift=-2.30cm,font=\color{black}]  {$t$};
  \node[above of= img1, node distance=0cm, yshift=-2.80cm,font=\color{black}]  {(a)};
  \node[left of= img1, node distance=0cm, rotate=90, anchor=center,yshift=2.60cm,font=\color{black}] { $x_t$};
\end{tikzpicture}\columnbreak
\hspace{0.15\textwidth}
\begin{tikzpicture}
  \node (img1)  {\includegraphics[width=0.274\textwidth]{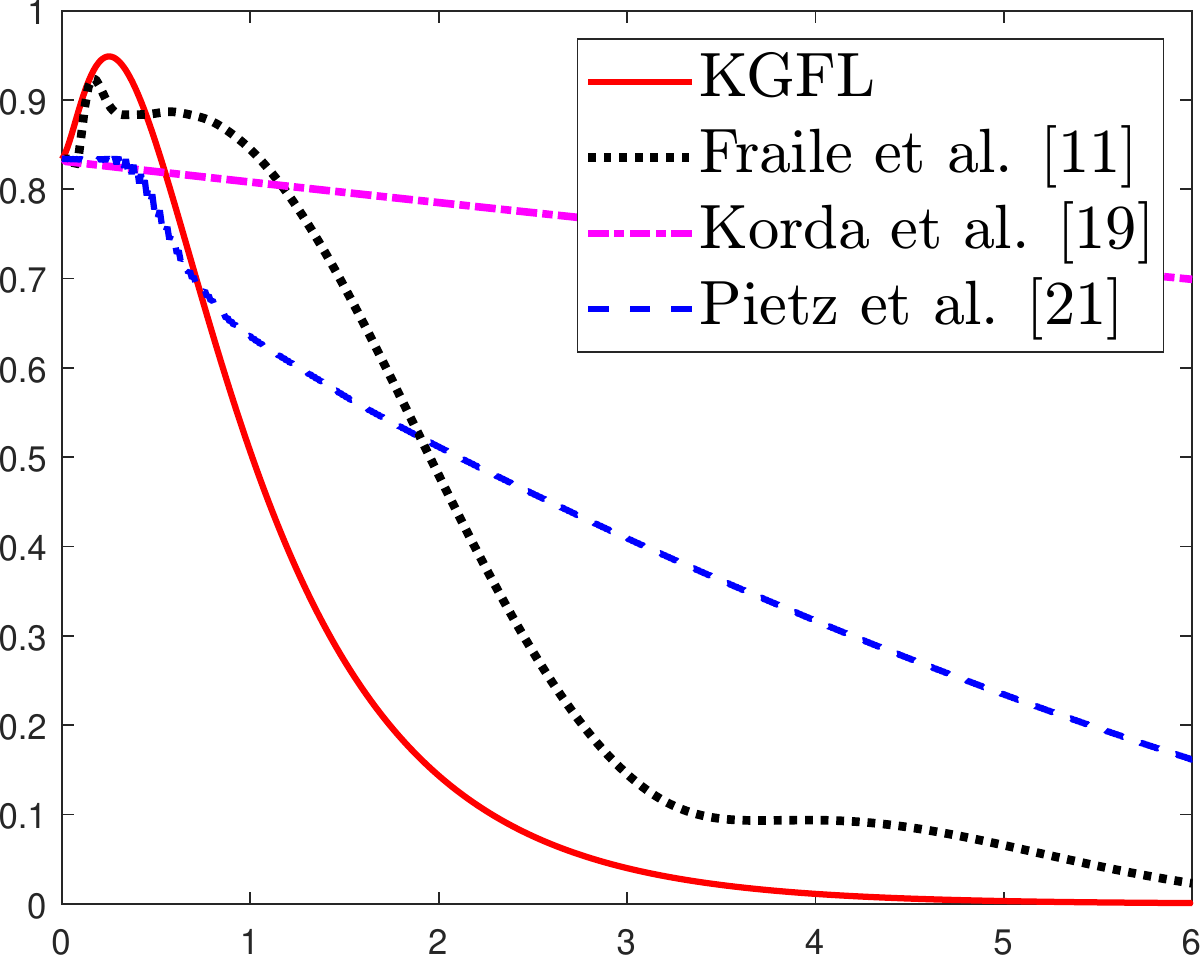}};
  \node[below of= img1, node distance=0cm, yshift=-2.3cm,font=\color{black}]  {$t$};
  \node[above of= img1, node distance=0cm, yshift=-2.8cm,font=\color{black}]  {(b)};
  \node[left of= img1, node distance=0cm, rotate=90, anchor=center,yshift=2.760cm,font=\color{black}] { $||x_t||_2$};
\end{tikzpicture}\columnbreak
\hspace{0.15\textwidth}
\begin{tikzpicture}
  \node (img1)  {\includegraphics[width=0.25\textwidth]{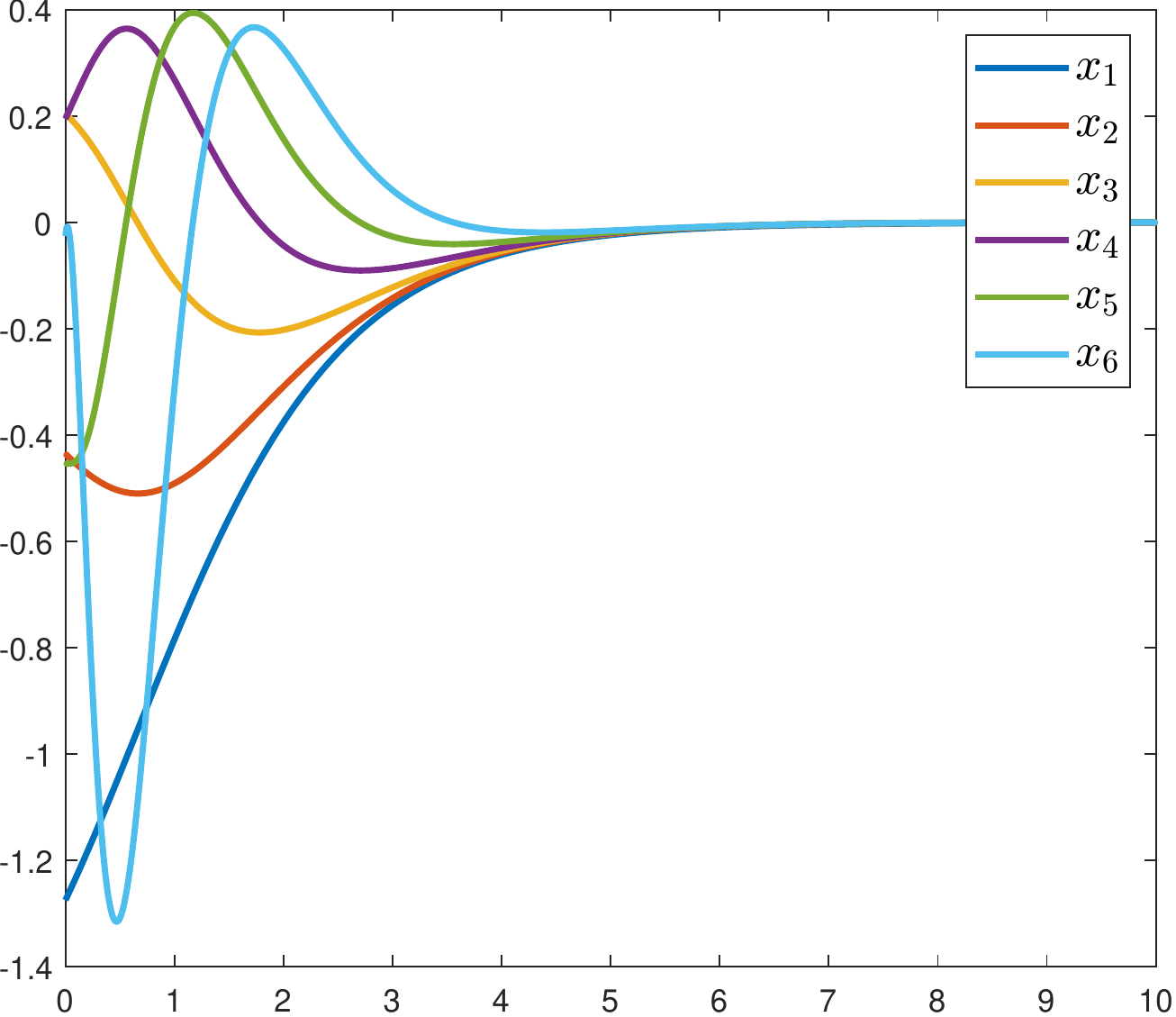}};
  \node[below of= img1, node distance=0cm, yshift=-2.3cm,font=\color{black}]  {$t$};
  \node[above of= img1, node distance=0cm, yshift=-2.8cm,font=\color{black}]  {(c)};
  \node[left of= img1, node distance=0cm, rotate=90, anchor=center,yshift=2.50cm,font=\color{black}] { $x_t$};
\end{tikzpicture}
\end{multicols}
\vspace*{-0.4cm}
\caption{\color{black}(a) Full state feedback linearization of the Van der Pol oscillator \eqref{eqn:vanderpol}. It
 can be seen that the data-driven algorithm(solid lines) performs almost as
 good as model-based feedback linearization (dashed lines) and much better than a simple linearization with no control transformation (dotted lines) \cite{MK-IM:2018}. The data-driven
 algorithm has slightly higher overshoot and settling time than the model-based
 method due to imperfections in the learned transformations. (b) Comparison of the proposed algorithm, \cite{SP-SEO-CWR:2020}, \cite{LF-MM-PT:2020} and \cite{MK-IM:2018}.
 (c) Full state feedback linearization for the complex 3 dimensional system \eqref{eqn:complex-sys}. The dictionary with Hermite polynomials of order 3 is augmented with $\sin(x)$. It is evident that the proposed algorithm stabilizes even a system of the third order. \label
 {fig:vanderpol}}
\end{figure*}

\begin{figure}
\begin{multicols}{2}
\begin{tikzpicture}
  \node (img1)  {\includegraphics[width=0.21\textwidth]{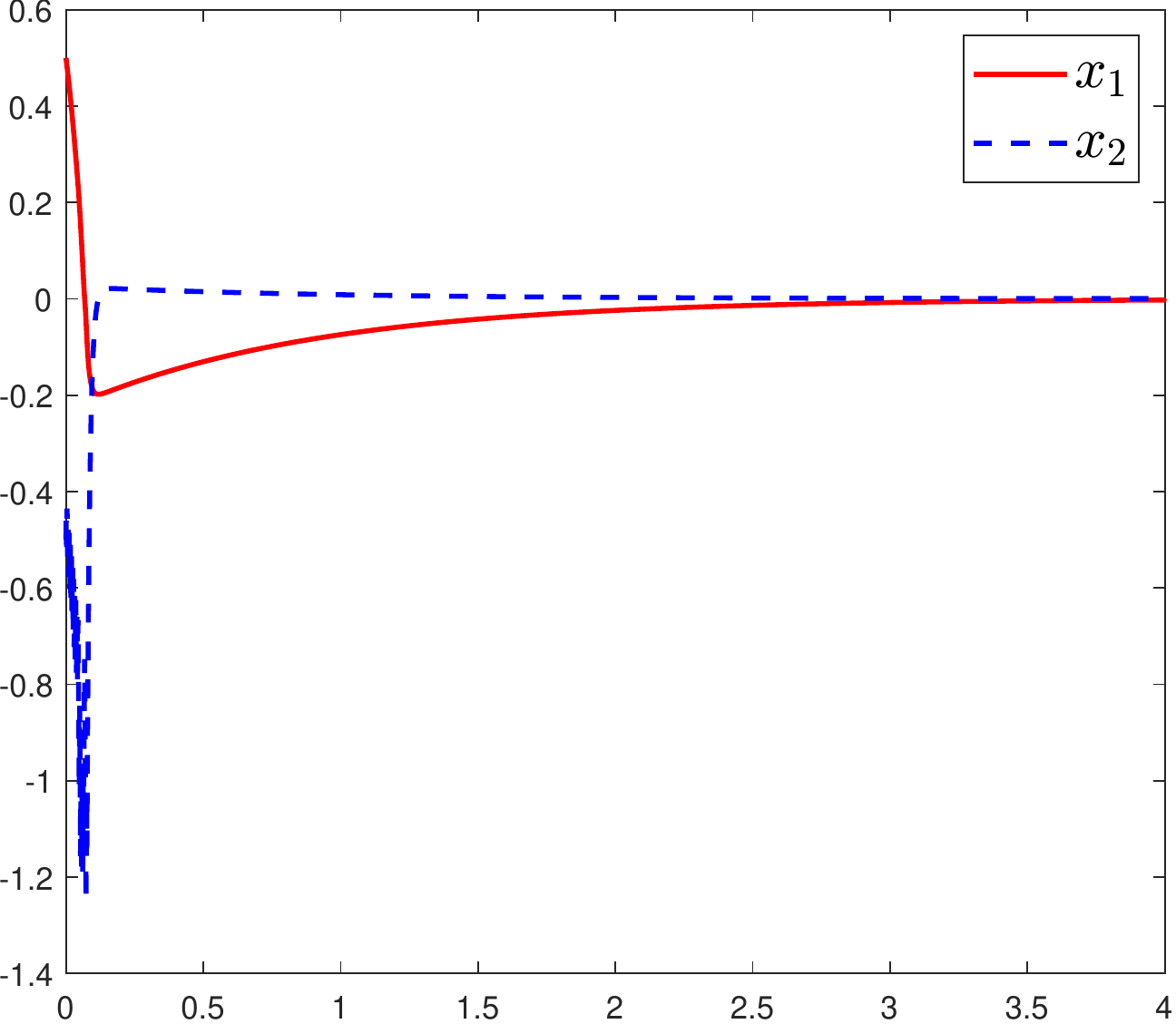}};
  \node[below of= img1, node distance=0cm,xshift = 0.2cm, yshift=-2.1cm,font=\color{black}]  {\footnotesize $t$};
  \node[above of= img1, node distance=0cm, yshift=-2.75cm,font=\color{black}]  {(a)};
  \node[left of= img1, node distance=0cm, rotate=90, anchor=center,yshift=2.1cm,font=\color{black}] { \small $x_t$};
\end{tikzpicture}\columnbreak
\hspace{0.1\textwidth}
\begin{tikzpicture}
  \node (img1)  {\includegraphics[width=0.21\textwidth]{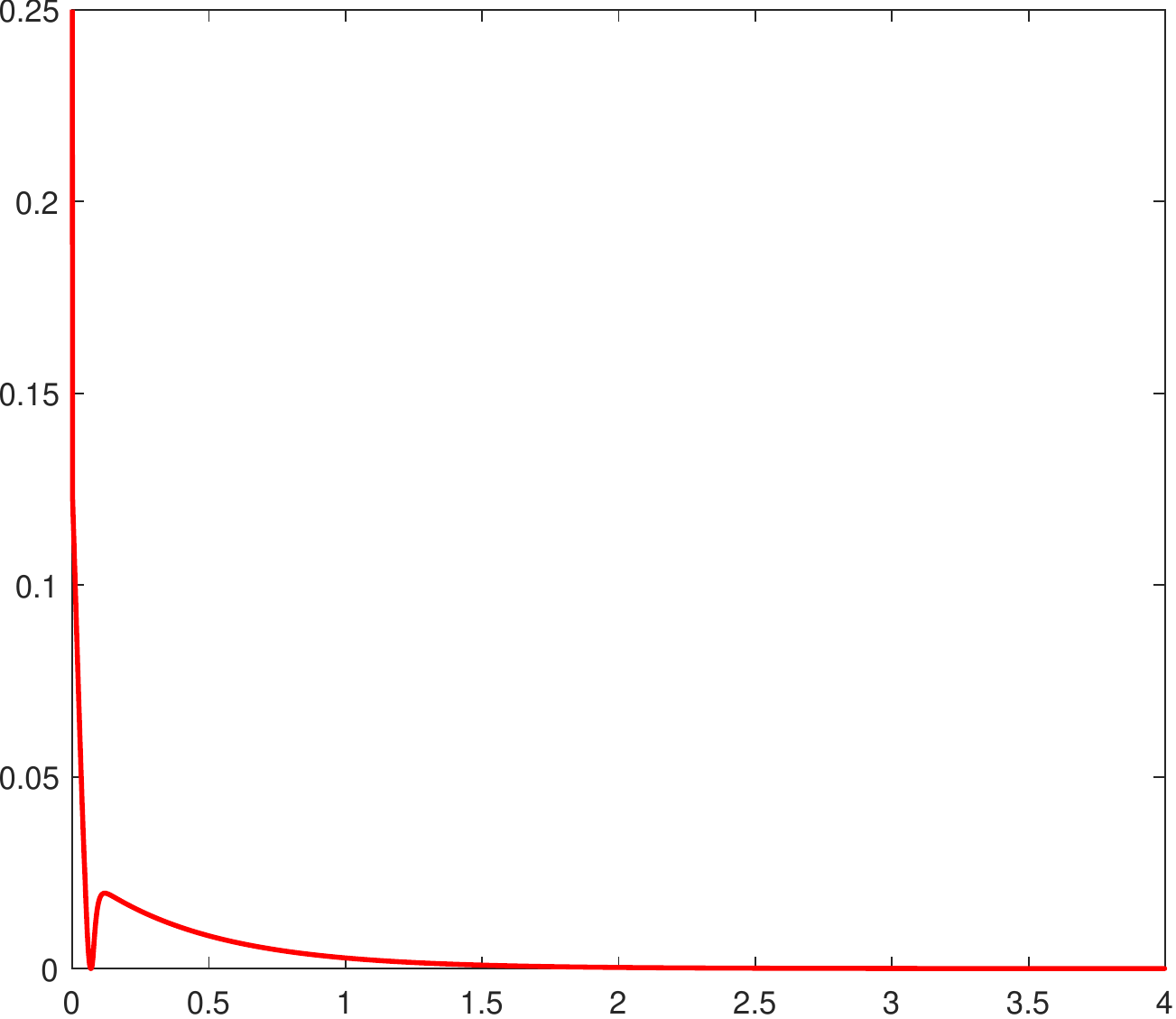}};
  \node[below of= img1, node distance=0cm, yshift=-2.1cm,font=\color{black}]  {$t$};
  \node[above of= img1, node distance=0cm, yshift=-2.75cm,font=\color{black}]  {(b)};
  \node[left of= img1, node distance=0cm, rotate=90, anchor=center,yshift=2.1cm,font=\color{black}] { \small $y_t$};
\end{tikzpicture}
\end{multicols}
\vspace*{-0.4cm}
\caption{{\color{black}(a) State trajectory for input-output feedback linearization of the modified Van der Pol oscillator \eqref{eqn:vanderpol}. The states are plotted on the y-axis whereas the time is plotted on the x-axis. (b) Trajectory of the output $y = 0.5 x_1^2$. The output is plotted on the y-axis. It is evident that both the output and the states are stabilized as feedback linearization of the output does not create any zero dynamics.}} \label{fig:iofl}
\end{figure}
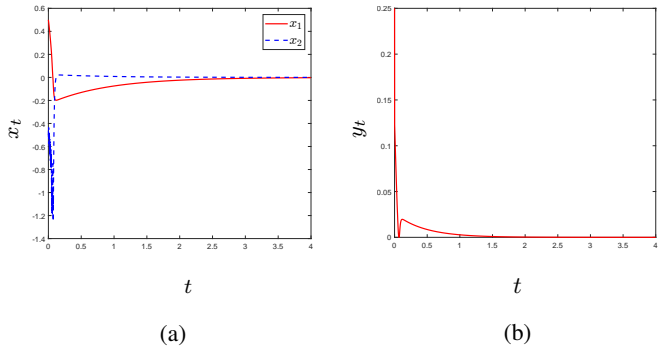

The Van der Pol oscillator is full-state feedback linearizable as
the system is both controllable and integrable. Particularly, choosing $h(x) =
x_{1}$, $\alpha(x) = x_{1} - 0.5(1-x_{1}^2)x_{2}$, and {\color{black}$\beta(x) = (1-x_2^2)^{-1}$} fully
linearizes the system in its normal form.
In Figure~\ref{fig:vanderpol}, the dashed lines represent model-based feedback
linearization whereas the solid lines represent the learned linearization
transformation using the proposed algorithm. In the figure, $x^m$ represents the
model-based states. In Figure~\ref{fig:vanderpol}(a), we choose a dictionary of
Hermite polynomials of order 2 and their mutual Kronecker products. We use KGFL to perform full-state feedback linearization. It is clear
from the figure that the data-driven algorithm is able to learn a
transformation and stabilize the system almost as well as the model-based
linearization. {\color{black} We also compare our proposed algorithm
 against the algorithms in \cite{MK-IM:2018}, \cite{SP-SEO-CWR:2020}, and \cite
{LF-MM-PT:2020} in Figure~\ref{fig:vanderpol}(b). In \cite{MK-IM:2018}, a linear predictor of the nonlinear system is constructed without any control transformations. Hence it cannot achieve exact linearization of the nonlinear system, especially when the control affects the system nonlinearly. The data-driven algorithm
proposed in \cite{LF-MM-PT:2020} is motivated by an intelligent PID controller
that makes use of sampled measurements in an online fashion. The algorithm
proposed in \cite{SP-SEO-CWR:2020} bilinearizes the system using offline data
and performs model predictive control. To maintain fairness in comparison, we
used the same initial conditions and the same offline data for the proposed
algorithm and \cite{SP-SEO-CWR:2020}. It
is evident from simulations that the proposed algorithm stabilizes faster than
the compared algorithms. The online algorithm in \cite{LF-MM-PT:2020} uses past
measurements to compute piece-wise constant inputs, whereas the algorithm
in \cite{SP-SEO-CWR:2020} performs MPC by bilinearizing the system by
performing state transformations. Our proposed algorithm linearizes the system
using both state and control transformations and applies a state-feedback
approach for pole placement which makes our algorithm effective. 
Furthermore, Figure~\ref{fig:vanderpol}(c) shows that the data-driven feedback linearization algorithm performs well even for the high dimensional system~\eqref{eqn:complex-sys} with inputs entering nonlinearly.

\subsubsection{Output feedback linearization}
We now demonstrate data-driven input-output feedback linearization on the 
Van der Pol oscillator using the single-step method described in Theorem~\ref{thm:state_io-fl}. We choose the output $y = h(x) = 0.5x_{1}^2$. We make the observation that the system has a relative degree $r = 2$ for the selected output. For the data-driven setting, we choose a
dictionaries that contains the Hermite polynomials of the output and its
derivatives up to the third degree. We learn the control transformation in
problem~\eqref{prob:optim-iofl}. In Figure~\ref{fig:iofl}(b), it is observed
that the proposed data-driven algorithm stabilizes the output at 0. Further,
the feedback linearization of the selected output induces no zero
dynamics, hence the overall system is also stable, as it can be observed in
figure~\ref{fig:iofl}(a). 
}

\subsection{Effect of richness of dictionary}

Here, we investigate the effect of the richness of the dictionary on input-output
feedback linearization for the single-step method. We choose the modified Van der Pol oscillator
introduced in equation~\eqref{eqn:vanderpol} and the output function $h
(x) = x_{1}$. Hence the state transformation becomes $z = \left[x_{1} \ x_
{2}\right]^\top$, and we only learn the control transformations $\hat{\alpha}$ and
$\hat{\beta}$. Note that we only use input-output data as described in section~\ref{subsect:fl}. We compare the difference in trajectories between that of model-based
feedback linearization and the proposed data-driven algorithm. We look at the
loss $Q_T$ defined as the sum of the norm difference between the states at every time
instance over a horizon $T$. That is, $Q_T = \sum\limits_t^T
||x^m_t - x_t||^2_2$. We choose $T = 10$ and use a data set of size $N = 300$.
We compare the loss in Figure~\ref{fig:comparison}(a) over 40 separate
experiments for each size of the dictionary. Figure~\ref{fig:comparison} is a
boxplot where the average is denoted by the red line and the $75^{th}$ and $25^
{th}$ quartiles are represented by upper and lower edges of the blue box
respectively. The red stars represent the outliers. We vary the richness of the
dictionary by varying the order of the Hermite polynomials included for learning
the control transformation. We see from Figure~\ref{fig:comparison}(a) that, as
the model complexity increases, the quantity $Q_T$ initially decreases but
starts to increase again after reaching a minimum. This can be attributed to
overfitting\cite{AMR-HMS:2020}. 


\subsection{Effect of dataset size}
\begin{figure*}
\begin{multicols}{3}
\begin{tikzpicture}
  \node (img1)  {\includegraphics[width=0.263\textwidth]{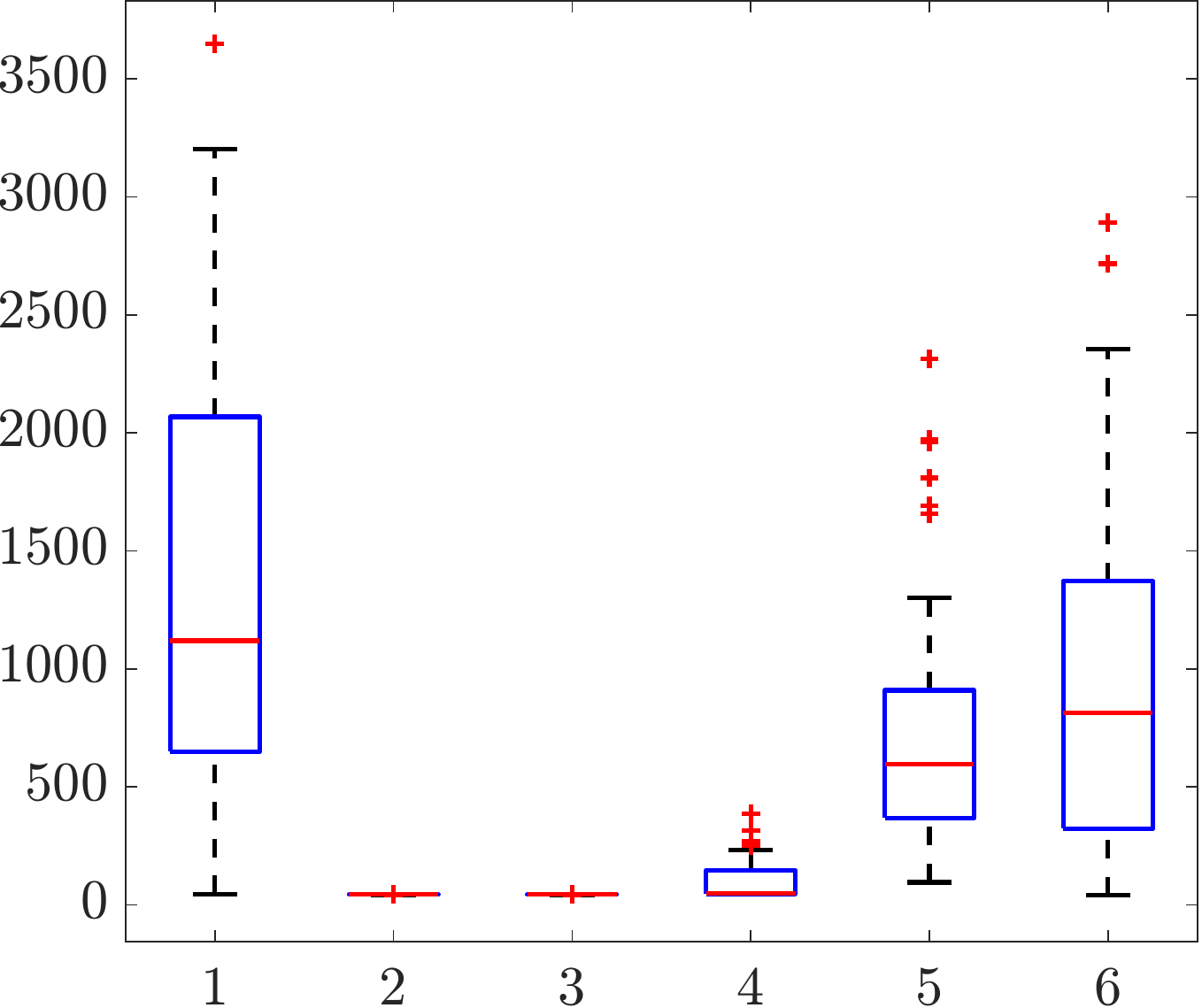}};
  \node[below of= img1, node distance=0cm,xshift = 0.2cm, yshift=-2.3cm,font=\color{black}]  { Order of Hermite polynomial};
  \node[above of= img1, node distance=0cm, yshift=-2.75cm,font=\color{black}]  {(a)};
  \node[left of= img1, node distance=0cm, rotate=90, anchor=center,yshift=2.6cm,font=\color{black}] { \small $\sum\limits_t ||x^m_t - x_t||^2_2$};
\end{tikzpicture}\columnbreak
\hspace{0.1\textwidth}
\begin{tikzpicture}
  \node (img1)  {\includegraphics[width=0.25\textwidth]{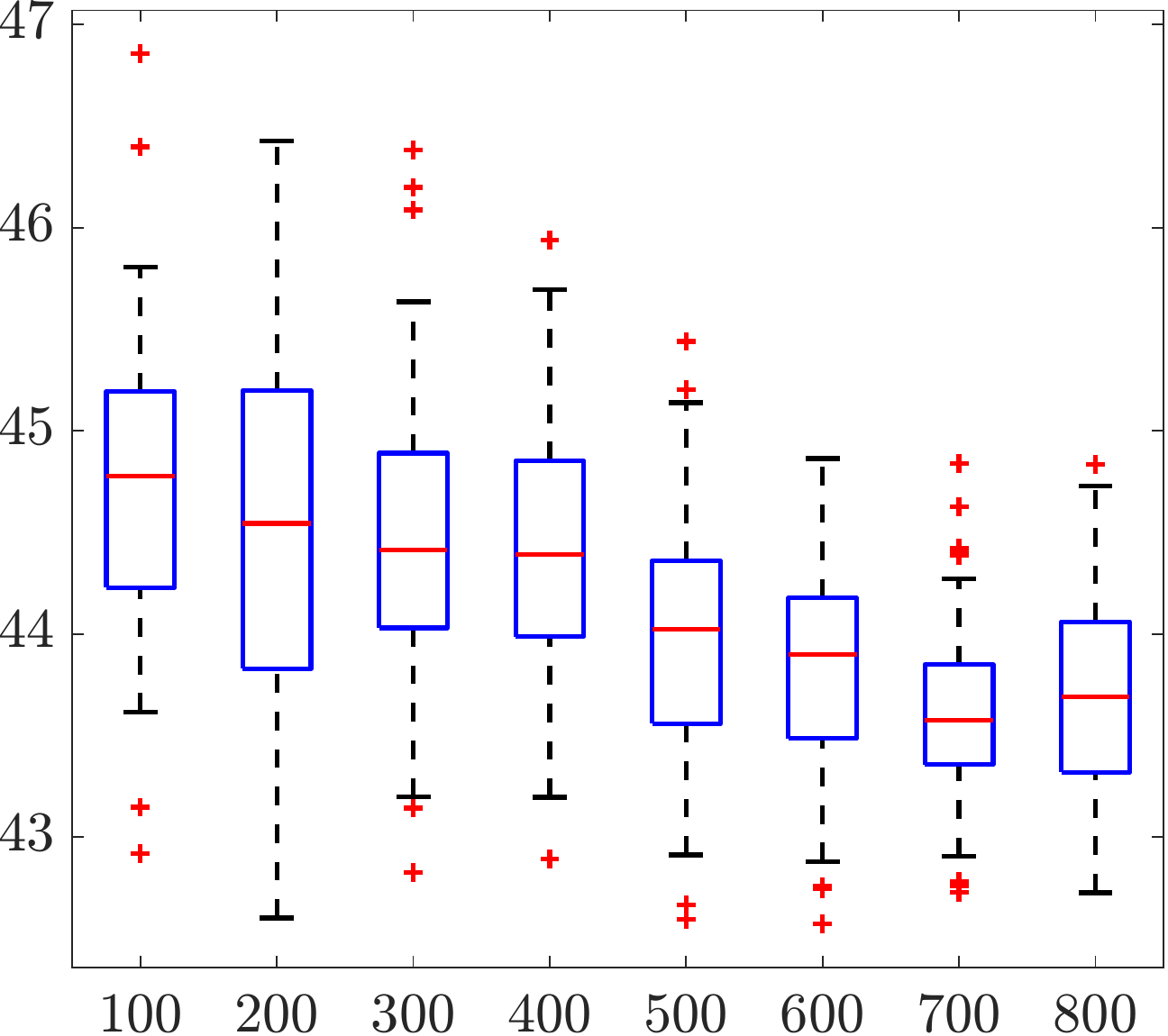}};
  \node[below of= img1, node distance=0cm, yshift=-2.3cm,font=\color{black}]  {Data size $N$};
  \node[above of= img1, node distance=0cm, yshift=-2.75cm,font=\color{black}]  {(b)};
  \node[left of= img1, node distance=0cm, rotate=90, anchor=center,yshift=2.5cm,font=\color{black}] { \small $\sum\limits_t ||x^m_t - x_t||^2_2$};
\end{tikzpicture}\columnbreak
\hspace{0.1\textwidth}
\begin{tikzpicture}
  \node (img1)  {\includegraphics[width=0.25\textwidth]{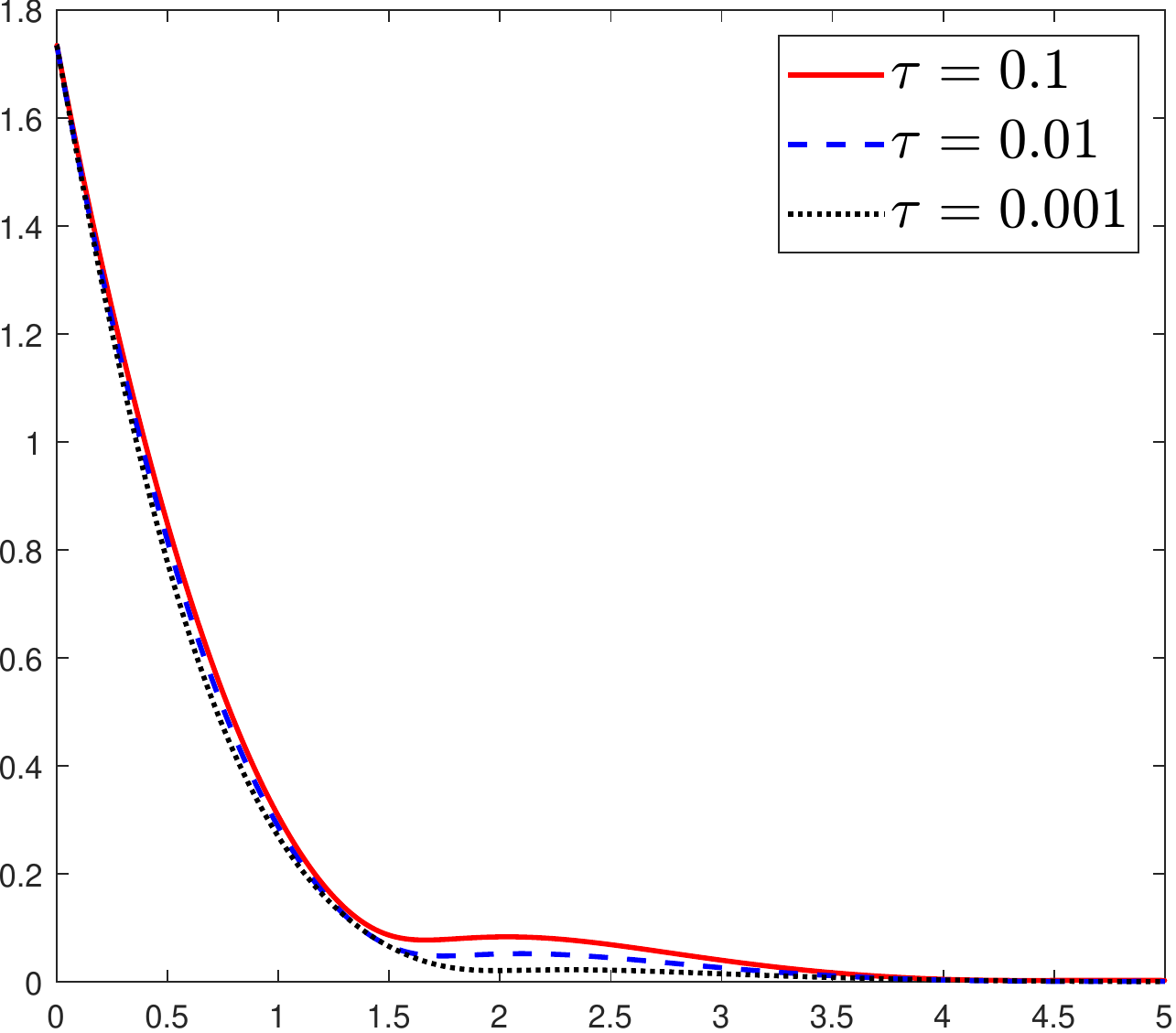}};
  \node[below of= img1, node distance=0cm, yshift=-2.3cm,font=\color{black}]  {$t$};
  \node[above of= img1, node distance=0cm, yshift=-2.75cm,font=\color{black}]  {(c)};
  \node[left of= img1, node distance=0cm, rotate=90, anchor=center,yshift=2.5cm,font=\color{black}] { \small $||x_t||$};
\end{tikzpicture}
\end{multicols}
\vspace*{-0.4cm}
\caption{(a) Richness of the dictionary vs. Loss. The deviation of the
 data-driven algorithm from model-based feedback linearization is considered as
 loss which is plotted on the y-axis. The orders of Hermite polynomials used in
 the dictionaries are plotted on the x-axis. As the dictionaries get richer, the
 loss initially decreases until polynomials of degree 3. However, as we
 consider richer dictionaries beyond order 4 Hermite polynomials, the model
 overfits the data and we incur high losses. (b) Data size vs. loss. The loss
 is defined similarly to the previous comparison. {\color{black}As the size of the
 data set increases, the loss uniformly decreases for a dictionary with Hermite
 polynomials up to the second degree. (c) The comparison of our algorithm KGFL for different sampling intervals. It can be seen that the sampling rate doesn't have a significant effect on the performance of the algorithm for the Van der Pol oscillator \eqref{eqn:vanderpol}. }} \label{fig:comparison}
\end{figure*}
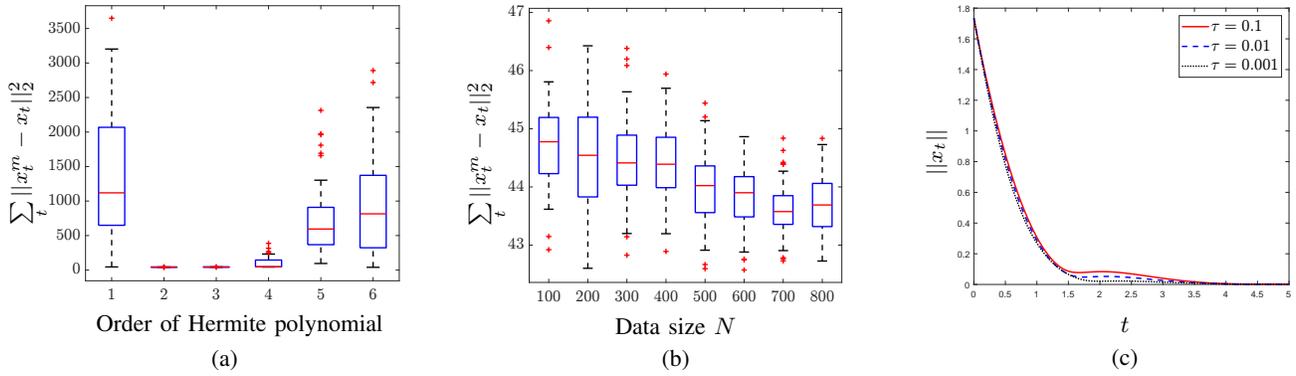

We study the effect of the size of the dataset on the loss $Q_T$ when the size of
the dictionary is fixed. Similar to the previous numerical comparisons, we perform
output feedback linearization with $y = h(x) = x_{1}$. {\color{black}The dictionaries for
control transformations contain Hermite polynomials up to the second degree.}
From Figure~\ref{fig:comparison}(b) we can see that as the data size increases,
the average loss uniformly decreases until 800 data points.


{\color{black}
\subsection{Effect of sampling interval}

We analyze the effect of sampling interval of the data on the performance of KGFL. For the Van der Pol oscillator considered in \eqref{eqn:vanderpol}, we sample the system at rates 0.1, 0.01, and 0.001 respectively. The results of the performance of these numerical experiments are presented in Figure \ref{fig:comparison}(c). It can be seen that the sampling rate does not have a significant effect on the performance of the algorithm. A fast sampling rate of 0.001 seems to perform only marginally better than sampling rates of 0.01 and 0.1. A complete analysis of the sampling rate is beyond the scope of the paper as it is dependent on the nature of the vector fields that describe the system. Typically, systems that evolve fast in the state space require faster sampling rates.  


}

\section{Conclusion} We establish a connection between the traditional
 model-based feedback linearization technique and the Koopman generator.
 Particularly, we show that here exists an observable and a state feedback
 control that renders the Koopman-generator finite-dimensional and nilpotent
 when the system is feedback linearizable. Using this connection, we develop an
 algorithm called KGFL to feedback linearize a control-affine system using
 experimental data. We demonstrate the algorithm numerically on complex
 dynamical systems and discuss tradeoffs related to the size of the
 dictionaries and the size of the dataset. We also show that it performs better
 than existing algorithms in the literature as KGFL exploits the feedback
 linearizable structure of the system. Directions of future research include
 the problem of choosing the right observables to obtain stable zero
 dynamics. 


\bibliographystyle{unsrt}
\bibliography{refs}

\begin{appendices}
\section{Proof of Theorem~\ref{thm:Koopman-nilpotent}}
\label{app:proof_koopman_nilpotent}
\textit{(i) Forward proof (Feedback linearization to Koopman generator): }Since $\left \lbrace \ad_{f}^{(k)} g \right \rbrace_{k=0}^{r-2}$ is involutive, Frobenius theorem guarantees that there exists $n-r+1$ functions $\left \lbrace h_i \right\rbrace_{i=1}^{n-r+1}$ that satisfy the following properties,
\begin{align*}
    L_g L_f^{k-1} h_i &= 0, \forall~k \in \lbrace 1, \ldots, r-1 \rbrace, \forall~i \in \lbrace 1, \ldots, n-r \rbrace, \\
    L_g L_f^{r-1} h_i &\neq 0.
\end{align*}
It then follows that for any $h \in \left\lbrace h_i \right\rbrace_{i=1}^{n-r+1}, k \in \lbrace 1, \ldots, r-1 \rbrace$ and any feedback $\alpha$,
we have $L_{f + g \alpha} L_f^{(k-1)} h = L_f^k h + \alpha L_g L_f^{k-1} h
= L_f^k h$. Now, let $\alpha(x) = - \left( L_g L_f^{r-1} h(x) \right)^{-1} L_f^r h(x)$.
It is important to note that $\alpha$ is well-defined at any $x$ since $ L_g L_f^{r-1} h(x) \neq 0$.
We then get $L_{f + g \alpha} L_f^{r-1} h
= L_f^r h + \alpha L_g L_f^{r-1} h = 0$.
Therefore, $L_{f + g\alpha}$ is locally nilpotent 
at~$h$, with index~$r$.

\textit{(ii) Converse proof (Koopman generator to feedback linearization):}
 First, we show that the subspace $\mathrm{span}\big( h, L_{f} h, \ldots, L_
 {f}^{r-1} h \big)$ is $L_{f + \alpha g}$ invariant.  Let $\mathcal
 {V} := \mathrm{span}\big( h, L_{f} h, \ldots, L_{f}^{r-1} h \big)$. For
 any $p \in \mathcal{V}$, we have $p = \sum_{k=1}^r c_k L_f^{(k-1)} h$, and it
 follows that:
\begin{align*}
    L_{f + g\alpha} p &= L_{f+g\alpha} \sum_{k=1}^r c_k L_{f}^{k-1} h  = \sum_{k=1}^r c_k L_{f + g\alpha} L_f^{k-1} h \\
    &= \sum_{k=1}^{r-1} c_k L_{f}^k h \in \mathcal{V}.
\end{align*}
Therefore, we get that $\mathcal{V}$ is $L_{f + g\alpha}$-invariant. Further, from Lemma 6.5 in \cite{JJES-WL-etal:1991} we have that functions $L_f^{i} h$ for $i=\{0,1,\dots,r-1\}$ are linearly independent. Therefore $\mathcal{V}$ is an $L_{f+g\alpha}$-cyclic basis.


We prove the converse of the theorem by induction. Let $\ell$ be an observable such that $L_{f+ g \alpha}$ is nilpotent at $\ell$ with index $r$. Further, let $\ell_k = L_{f}^{(k-1)} \ell$ for $k \in \lbrace 1, \ldots, r \rbrace$. Consider the functions $\ell_1$ and $\ell_2$. We have that $L_{f+g \alpha} \ell_1 = L_f \ell_1 + \alpha L_g \ell_1 = \ell_2$. However, we have that $\frac{\partial \ell_2}{\partial \alpha} = 0$ due to the stability of $L_{f+\alpha g}$-cyclic subspace, where $\frac{\partial}{\partial \alpha}$ denotes the Gateaux derivative with respect to $\alpha$. This is a consequence of the definition of subspace stability. This implies that $L_g \ell_1 = 0$. Similarly, considering $\frac{\partial \ell_3}{\partial \alpha} = 0$, we get that $L_g L_f \ell_1 = 0$. By doing this every basis function $\ell_i$, we obtain that $\mathrm{Ker}(L_g) = \mathrm{span} \left \lbrace \ell_1, \ell_2, \ldots, \ell_{r-1} \right \rbrace $ as $L_g \ell_i = 0$ for $i \in \{1, 2, \dots, r-1\}$. Therefore $\mathrm{Ker}(L_g)$ is $r-1$ dimensional. We know that $\mathcal{V}$ is $r$-dimensional and is $L_{f + g \alpha}$-invariant
subspace with a $L_{f + g \alpha}$-cyclic basis $\left \lbrace h_1, h_2, \ldots, h_r \right \rbrace$. Now, we can see that $\mathcal{V} \cap \mathrm{Ker}(L_g)$ is $(r-1)$-dimensional. This implies that the distribution $\{g,\ad_fg,\dots, \ad_f^{r-2}g\}$ is involutive. 

\section{Proof of Theorem~\ref{thm:state_io-fl}}
\label{app:proof_lemma_io-full}
\begin{proof}If the dictionaries $D$, $\theta$ and $\gamma$ are complete, we can express the feedback linearized system with data as follows. 
\begin{align*}
\dot{D}(X)K &= A D(X)K + B(Q^\top \Theta(Y) + R^\top \Gamma(Y) \otimes U).
\end{align*}
The least-squares solution in terms of $K$, $Q$ and $R$ can be obtained by simply solving for the first-order condition of optimality. For input-output feedback linearization, $\dot{D}(X)K$ and $D(X)K$ are replaced by $\mathbf{\dot{Z}}$ and $\mathbf{Z}$ respectively.
\end{proof}
\end{appendices}

\end{document}